\documentclass[amsfonts]{article}
\usepackage{pslatex}
\usepackage{amsmath}
\usepackage{amssymb}
\usepackage{amsfonts}
\usepackage{latexsym}
\begin{document}

\title{Kolmogorov - Sinaj entropy on MV-algebras}

\author{Beloslav Rie\v can
 \thanks{M. Bel University, Tajovsk\' eho 40, 97401
Bansk\' a Bystrica, Matematick\' y \' ustav SAV, \v Stef\' anikova
49, 81473 Bratislava,\ttfamily{riecan@fpv.umb.sk,
riecan@mat.savba.sk}}}
\date{}
\maketitle

%

 \begin{abstract}There is given a construction of the entropy
of a dynamical system on arbitrary MV-algebra $M$. If $M$ is the
MV-algebra of characteristic functions of a $\sigma$-algebra
(isomorphic to the $\sigma$-algebra), then the construction leads
to the Kolmogorov - Sinaj entropy. If $M$ is the MV-algebra
(tribe) of fuzzy sets, then the construction coincides with the
Mali\v ck\' y modification of the Kolmogorov - Sinaj entropy for
fuzzy sets ([6], [14], [15]).
\end{abstract}

\section{ Introduction}

If $(\Omega, \mathcal{S}, P)$ is a probability space and $T:\Omega
\rightarrow \Omega$ is  a measure preserving transformation
(i.e. $A \in \mathcal{S}$ implies $T^{-1}(A) \in \mathcal{S}$, and
$P(T^{-1}(A)) = P(A)$, the entropy is defined as follows. If
$\mathcal{A} = \{A_1,...,A_k\}$ is a measurable partition of
$\Omega$, then
$$
H(A) = \sum_i \varphi (P(A_i)),
$$
where $\varphi (x) = - x log x$, if $x > 0$, and $\varphi (0) =
0$. Denote by $\bigvee_{i=0}^{n-1}T^{-i}(\mathcal{A})$ the common
refinement of the partitions $\mathcal{A}, T^{-1}(\mathcal{A}) =
\{T^{-1}(A_1),...,T^{-1}(A_k)\}, ..., T^{-(n-1)}(\mathcal{A}) =
\{T^{-(n-1)}(A_1),...,T^{-(n-1)}(A_k)\}.$ Then there exists
\begin{equation}
\label{rov 1} h(\mathcal{A},T) = \lim_{n \to \infty} \frac 1 n
H(\bigvee_{i=0}^{n-1} T^{-i}(\mathcal{A}))
\end{equation}
and the Kolomogorov - Sinaj entropy is defined by the formula

 \begin{equation}h(T) = \sup \{ h(\mathcal{A},T); \mathcal{A} \ {\text{is a
measurable partition of}} \ \Omega \}.\label{rov 2}\end{equation}

Let us try to substitute partitions by fuzzy partitions
$\mathcal{A} = \{f_1,...,f_k\}$, where $f_j$ are non-negative,
measurable functions, $\sum_{i=1}^k f_i = 1_{\Omega}.$ Put \ \
$m(f_i) = \ \int_{\Omega} f_i $ dP. The common refinement of
$\mathcal{A}$ and $\mathcal{B} = \{ g_1,...,g_l\} $ is defined
 using products of functions: $\mathcal{A} \vee \mathcal{B} = \{f_i.g_j; i=1,...,k,
j=1,...,l\}$. Further $T^{-1} (\mathcal{A}) = \{ f_1 \circ T,...,
f_k \circ T\}.$ Then the entropy can be defined by the formula
$$
H(\mathcal{A}) = \sum_{i=1}^k \varphi (m(f_i)),
$$
and then by (1) and (2). Of course, in this case $h(T) = \infty$,
if e.g. the supremum is taken over the all measurable partitions.
Therefore P. Mali\v ck\'y suggested ([6], see also [14], [15]) the following modification:
instead of $H(\bigvee_{i=0}^{n-1} T^{-1}(\mathcal{A})$ to use

 \begin{equation}\label{rov 3}H_n(\mathcal{A}, T) = \inf \{ H(\mathcal{C}); \mathcal{C} \geq
\mathcal{A}, \mathcal{C} \geq T^{-1}(\mathcal{A}),..., \mathcal{C}
\geq T^{-(n-1)}(\mathcal{A})\}\end{equation}
 and then
to define the entropy by the formula

 \begin{equation}\label{rov 4}h(\mathcal{A},T) = \lim_{n \to \infty} \frac 1 n
H_n(\mathcal{A},T).\end{equation}

A natural algebraization of a tribe of functions is an MV-algebra. On a
special type of MV-algebras - so-called MV-algebras with product
([8],[13]) there has been realized the construction of entropy in [11],[12].
In the paper, following the Mali\v ck\' y construction (3), (4), we
construct entropy on any MV-algebra. Secondly, we need only additivity of
considered states instead of usually demended $\sigma$-additivity.

In Section 2 we present necessary informations about MV-algebras
and some auxi\-liary results concerning partitions of unity. Then
in Section 3 we prove the existence of entropy. Section 4 registers basic facts
about isomorphism and entropy and Section 5 contains a counting formula.

\section{ MV-algebras}

An MV-algebra $M = (M,0,1,\neg, \oplus, \odot)$ is a system where
$\oplus$ is associative and commu\-tative with neutral element
$0$, and, in addition, $\neg 0 = 1, \neg 1 = 0, x \oplus 1 = 1, x
\odot y = \neg (\neg x \oplus \neg y),$ and $y \oplus \neg (y
\oplus \neg x) = x \oplus \neg (x \oplus \neg y)$ for all $x,y \in
M$.

MV-algebras stand to the infinite-valued calculus of Lukasiewicz as
boolean algebras stand to classical two-valued calculus.

An example of an MV-algebra is the real unit interval [0,1] equipped with the
ope\-rations
$$
\neg x = 1 - x, x\oplus y = \min (1,x + y), x \odot y = \max (0,x + y - 1)
$$
It is interesting that any MV-algebra has a similar structure. Let $G$ be
a lattice-ordered Abelian group (shortly $l$-group). Let $u \in G$ be a
strong unit of $G$, i.e. for all $g \in G$ there exists and integer $n
\geq 1$ such that $nu \geq g$. Let $\Gamma (G,u)$ be the unit interval
$[0,u] = \{h \in G; 0 \leq h \leq u \} $ equipped with the operations
$$
\neg g = u - g, g \oplus h = u \wedge (g + h), g \odot h = 0 \vee (g + h -
u).
$$
Then $([0,u], 0, u, \neg, \oplus, \odot)$ is an MV-algebra and by the
Mundici theorem ([9]), up to isomorphism, every MV-algebra $M$ can be
identified with the unit interval of a unique $l$-group $G$ with strong
unit, $M = \Gamma (G,u)$.

A partition of unity $u$ in $M$ is an $k$-tuple $\mathcal{A} =
(a_1,...,a_k)$ of elements of $M$ such that
$$
a_1 + ... + a_k = u.
$$

If $\mathcal{A} = \{a_1,...,a_n\}$, and $\mathcal{B} =
\{c_1,....,c_k\}$ are partitions of unity, then their common
refinement is any matrix  $\mathcal{S}=\{c_{ij}; i=1,...,k,
j=1,...l\} $ of elements of $M$ such that
$$
a_i = \sum_{j=1}^{l} c_{ij}, i = 1,...,k, b_j = \sum_{i=1}^{k}
c_{ij}, j = 1,...,l
$$
\bigskip
Although in view of the Mundici theorem the next result is known ([5] Prop. 2.2 (c), see also [7]), we mention the proof.
\bigskip

\bf Lemma 1. \it To any partitions $\mathcal{A}, \mathcal{B}$
there exists their common refinement.

\rm Proof. Let $\mathcal{A} = \{a_1,a_2,...,a_n\}, \mathcal{B} =
\{b_1,b_2,..., b_m\}$. By the well known Riesz theorem, to any
$a,b,c$ such that $0 \leq a \leq b + c, b \geq 0, c \geq 0$ there
exist $d,e$ such that $0 \leq d \leq b, 0 \leq e \leq c,$ and $a =
d + e$. (Indeed, put $d = a \wedge b, e = a - a \wedge b.)$ This
pro\-perty can be generalized by induction for a finite number of
elements. Therefore, since $a_1 \leq b_1 + b_2 + ... + b_m = u$,
there exist $c_{11}, c_{12},...,c_{1m}$ such that
$$
c_{1i} \leq b_i (i=1,....,m),
$$
and
$$
a_1 = c_{11} + c_{12} + ... + c_{1m}.
$$
Since $a_2 \leq u - a_1 = (b_1 - c_{11}) + (b_2 - c_{12}) + ... + (b_m -
c_{1m})$, there exist $c_{2i} (i=1,2,...,m)$ such that
$$
c_{2i} \leq b_i - c_{1i} (i=1,2,...,m),
$$
and
$$
a_2 = c_{21} + c_{22} + ... + c_{2m}.
$$
Continuing the process we construct $c_{3i} (i=1,2,...,m)$ such that
$$
c_{3i} \leq b_i - c_{1i} - c_{2i} (i=1,2,...,m),
$$
and
$$
a_3 = c_{31} + c_{32} + ... + c_{3m}.
$$
After constructing $c_{n-1,1}, c_{n-1,2},... c_{n-1,m}$ with $c_{n-1,i}
\leq b_i - \sum_{j=1}^{n-2}c_{ji}, a_{n-1} = \sum_{j=1}^m c_{n-1,j}$, we
put
$$
c_{ni} = b_i - \sum_{j=1}^{n-1} c_{ji} (i=1,...,m)
$$
Then evidently
$$
b_i = \sum_{j=1}^n c_{ji} (i=1,...,m),
$$
and
$$
c_{n1} + c_{n2} + ... + c_{nm} = b_1 + ... + b_m - \sum_{i=1}^m
\sum_{j=1}^{n-1} c_{ji}
$$
$$
 = u - \sum_{j=1}^{n-1} \sum_{i=1}^m c_{ji} = u - \sum_{j=1}^{n-1} a_j
= a_n.
$$
\bigskip
\bf Remark. \rm Common refinements of partitions are not defined
uniquely. E.g. let $M = [0,1], \mathcal{A} = \{0.5, 0.5\},
\mathcal{B} = \{0.4, 0.6\}.$ Then any matrix
$$
t - 0.1;  0.5 - t
$$
$$
0.6 - t; \ \ t
$$
$t \in [0.1, 0.5]$ represents a common refinement of $\mathcal{A},
\mathcal{B}$.

\section{ Entropy of dynamical systems}

\bf Definition 1. \it By a dynamical system on an MV-algebra we understand
a couple of mappings $m:M \rightarrow [0,1], \tau:M \rightarrow M$
satisfying the following conditions:
\begin{itemize}

\item[(i)]$
if a = b + c, then \ m(a) = m(b) + m(c), \tau (a) = \tau (b) +
\tau (c)$

\item[(ii)]$
\tau (u) = u, m(u) = 1$

\item[(iii)] $m(\tau (a)) = m(a), a \in M
$
\end{itemize}
\bf Definition 2. \it If $\mathcal{A} = \{ a_1,...,a_n\}$ is a
partition of unity, then its entropy is defined by the formula
$$
H(\mathcal{A}) = \sum_{i=1}^n \varphi (m(a_i)),
$$
where $\varphi (x) = - x log x$, if $ x > 0, \varphi (0) = 0$.

\bf Definition 3. \it If $\mathcal{A} = \{ a_1,...,a_n\},
\mathcal{B} = \{ b_1,...,b_k\}$ are two partitions of unity and
$\mathcal{C} = \{ c_{ij}; i=1,...,n, j = 1,...,m\}$ is a
common refinement of $\mathcal{A}$ and $\mathcal{B}$, then we define
\bigskip
$H_{\mathcal{C}}(\mathcal{A} \vert \mathcal{B})
 = \sum_{i=1}^n \sum_{j=1}^m m(b_j)\varphi(\frac {m(c_{ij})} {m(b_j)} ).$
\bigskip

\bf Lemma 2. \it $ H_{\mathcal{C}}(\mathcal{A} \vert \mathcal{B})
\leq H(\mathcal{A})$.

\rm Proof. Fix $i$, and put $\alpha_j = m(b_j), x_j = \frac {m(c_{ij})}
{m(b_j)}, j=1,...,k.$ Then $\sum_{j=1}^m \alpha_j = m(u) = 1$
$$
\sum_{j=1}^k \alpha_j x_j = \sum_{j=1}^k m(c_{ij}) = m(\sum_{j=1}^k c_{ij})
= m(a_i).
$$
Since $\varphi$ is concave, we have
$$
\sum_{j=1}^k \alpha_j \varphi (x_j) \leq \varphi (\sum_{j=1}^k \alpha_j
x_j),
$$
hence
$$
H_{\mathcal{C}}(\mathcal{A} \vert \mathcal{B}) = \sum_{i=1}^n
\sum_{j=1}^k m(b_j) \varphi ( \frac {m(c_{ij})} {m(b_j)}) =
\sum_{j=1}^n \sum_{j=1}^k \alpha_j \varphi (x_j) \leq \sum_{i=1}^n
\varphi (\sum_{j=1}^k \alpha_j x_j) = \sum_{i=1}^n \varphi
(m(a_i)) = H(\mathcal{A}).
$$

\bf Lemma 3. \it $H(\mathcal{C}) = H(\mathcal{A}) +
H_{\mathcal{C}}(\mathcal{B} \vert \mathcal{A})$ for any common
refinement $\mathcal{C}$ of $\mathcal{A}$ and $\mathcal{B}$.

\rm Proof. It follows by the additivity of logarithms.
$$
H(\mathcal{C}) = \sum_{i=1}^n \sum_{j=1}^k \varphi (c_{ij}) =
\sum_{i=1}^n \sum_{j=1}^k \varphi (m(a_i) \frac {m(c_{ij})}
{m(a_i)})
$$
$$
= - \sum_i \sum_j m(c_{ij}) log m(a_i) - \sum_i \sum_j m(c_{ij} \log \frac
{m(c_{ij})} {m(a_i)}
$$
$$
= - \sum_i \log m(a_i) \sum_j m(c_{ij}) + \sum_i \sum_j m(a_i)
\varphi (\frac {m(c_{ij})} {m(a_j)})
$$
$$
= - \sum_i m(a_i) \log m(a_i) + H_{\mathcal{C}} (\mathcal{B} \vert
\mathcal{A})
$$
$$
= H(\mathcal{A}) + H_{\mathcal{C}}(\mathcal{B} \vert \mathcal{A})
$$

\bf Corollary. \it For any common refinement $\mathcal{C}$ of
$\mathcal{A}$ and $\mathcal{B}$ there holds $H(\mathcal{C}) \leq
H(\mathcal{A}) + H(\mathcal{B}).$

\bf Lemma 4. \it For any partition $\mathcal{A} =
\{a_1,....,a_n\}$ put $\tau (\mathcal {(A)} =
\{\tau(a_1),...,\tau(a_n)\}$. Then $\tau(\mathcal{A})$ is a
partition of unity, too. Moreover, $H(\tau(\mathcal{A})) =
H(\mathcal{A})$.

\rm Proof. By (i) and (iii) of the definition of dynamical system we have
$$
\tau(a_1) + \tau(a_2) + ... + \tau (a_n) = \tau(a_1 + ... + a_n) = \tau(u) =
u.
$$
Moreover, by (iii)
$$
H(\tau(\mathcal{A})) = \sum_{i=1}^n \varphi (m(\tau(a_i))) =
\sum_{i=1}^n \varphi (m(a_i)) = H(\mathcal{A}).
$$

\bf Definition 3. \it For any partition $\mathcal{A}$ of unity and
any positive integer $n$ we define
\bigskip
$H_n(\mathcal{A}, \tau) = \inf \{ H(\mathcal{C}); \mathcal{C}$ is a common refinement of $\mathcal{A}, \tau(\mathcal{A}),..., \tau^{(n-1)}(\mathcal{A})\}$.
\bigskip

\bf Theorem 1. \it There exists $ \lim_{n \to \infty} \frac 1 n
H_n(\mathcal{A},T).$

\rm Proof. It suffices to prove that $H_{n+m}(\mathcal{A},\tau)
\leq H_n(\mathcal{A}, \tau) + H_m(\mathcal{A}, \tau)$.We shall use the following notation:  $\mathcal{C} \in \mathcal{A} \vee \mathcal{B}$, if $\mathcal{C}$ is a common refinement of $\mathcal{A}$ and $\mathcal{B}$. Let
$\mathcal{C}, \mathcal{D}$ be partitions of unity such that

\begin{equation}\label{rov 5}\mathcal{C} \in \mathcal{A} \vee \tau(\mathcal{A}) \vee
... \vee \tau^{(n-1)}(\mathcal{A}),\end{equation}

\begin{equation}\label{rov 6}\mathcal{D} \in \mathcal{A} \vee \tau(\mathcal{A}) \vee
... \vee \tau^{(m-1)}(\mathcal{A}).\end{equation}

By (6) we obtain

\begin{equation}\label{rov 7}\tau^n (\mathcal{D}) \in \tau^n(\mathcal{A}) \vee
\tau^{n+1}(\mathcal{A}) \vee ... \vee
\tau^{n+m-1}(\mathcal{A}).\end{equation}

Finally consider any partition $\mathcal{E}$ such that

 \begin{equation}\label{rov 8}\mathcal{E} \in \mathcal{C} \vee
\tau_n(\mathcal{D}).\end{equation}

By (8), (5) and (7) we obtain
$$
\mathcal{E} \in \mathcal{A} \vee \tau (\mathcal{A}) \vee
... \vee \tau^{n+m-1}(\mathcal{A}),
$$
hence
$$
H_{n+m}(\mathcal{A}) \leq H(\mathcal{E}).
$$
Of course, by (8), Corollary and Lemma 4
$$
H(\mathcal{E}) \leq H(\mathcal{C}) + H(\tau_n(\mathcal{D}))
$$
$$ = H(\mathcal{C}) + H(\mathcal{D}),
$$
hence

\begin{equation}\label{rov 9}H_{n+m}(\mathcal{A}) \leq H(\mathcal{C}) + H(\mathcal{D}).\end{equation}

Fix for a moment $\mathcal{D}$. Since (9) holds for any
$\mathcal{C} \in \mathcal{A} \vee
\tau(\mathcal{A})\vee ... \vee \tau^{(n-1)}(\mathcal{A})$, by the
definition of $H_n(\mathcal{A})$ we obtain
$$
H_{n+m}(\mathcal{A}) - H(\mathcal{D}) \leq H_n(\mathcal{A}),
$$
and by a similar argument
$$
H_{n+m}(\mathcal{A}) - H_n(\mathcal{A}) \leq H_m(\mathcal{A}).
$$
\bigskip

\bigskip

\section{ Isomorphism}

\bigskip

\bigskip

The aim of the Kolmogorov - Sinaj entropy was to distinguish non
-isomorphic dynamical systems. If $T \sim T'$ implies $h(T)= h(T')$, then
dynamical systems $T, T'$ with different entropies $h(T) \neq h(T')$
cannot be isomorphic.

\bigskip

\bf Definition 4. \it Entropy of an MV-algebra dynamical system
$(M,m,\tau)$ is defined by the formula
$$
h(\tau) = \sup \{ h(\mathcal{A},\tau); \mathcal{A} \ \text {is a
partition of unity}\}
$$
where
$$
h(\mathcal{A}, \tau) = \lim_{n \to \infty} \frac 1 n
H_n(\mathcal{A}, \tau).
$$

\bigskip

\bf Definition 5. \it Two MV-algebra dynamical systems
$(M_1,m_1,\tau_1), (M_2,m_2,\tau_2)$ are equivalent, if there
exists a mapping $\psi :M_1 \rightarrow M_2$ satisfying the
following conditions:
 \begin{itemize}
\item[(i)]  $ \psi \ \text {is a
bijection} $
\item[(ii)]
$ if \ a,b,c \in M_1, a = b + c,\text {then}\ \psi (a) = \psi (b)
+ \psi (c). $

\item [(iii)]$ \psi (u_1) = u_2 $
\item [(iv)]$ m_2(\psi(a)) = m_1(a) \text { for any}\ a \in M_1$
\item[(v)]
$ \tau_2(\psi(a)) = \psi(\tau_1(a)) \text { for any}\ a \in M_1 $
\end{itemize}
\bigskip

\bigskip

\bf Theorem 2. \it If $(M_1,m_1,\tau_1)$ and $(M_2,m_2,\tau_2)$ are
equivalent, then $h(\tau_1) = h(\tau_2)$.

\bigskip

\rm Proof. If $\mathcal{A}$ is any partition of $u_1$, then
$\psi(\mathcal{A})$ is a partition of $u_2$, and $H(\mathcal{A}) =
H(\psi(\mathcal{A}))$. Let $\varepsilon$ be an arbitrary positive
number. Choose a common refinement $\mathcal{C}$ of $\mathcal{A}, ...,
 \tau_1^{n-1}(\mathcal{A})$ such that
$$
H_n + \varepsilon > H(\mathcal{C}).
$$
Evidently $H(\mathcal{C}) \geq H_n(\psi(\mathcal{A}))$. Since
$$
H_n(\mathcal{A}) + \varepsilon \geq H_n(\psi(\mathcal{A}))
$$
holds for any $\varepsilon > 0$, we have
$$
H_n(\mathcal{A}) \geq H_n(\psi(\mathcal{A})),
$$
hence
$$
h_1(\tau_1,\mathcal{A}) = \lim_{n \to \infty} \frac 1 n
H_n(\mathcal{A}) \geq \lim_{n \to \infty} \frac 1 n
H_n(\psi(\mathcal{A})) = h_2(\tau_2, \psi(\mathcal{A})),
$$
$$
h_1(\tau_1) = \sup \{ h_1(\tau_1,\mathcal{A}); \mathcal{A}\} \geq
h_2(\tau_2, \psi(\mathcal{A})).
$$
Let $\mathcal{B}$ be any partition of $u_2$. Then $\mathcal{A} =
\psi^{-1}(\mathcal{B})$ is a partition of $u_1$ and
$\psi(\mathcal{A}) = \mathcal{B}.$ Therefore
$$
h_1(\tau_1) \geq h_2(\tau_2,\psi(\mathcal{A})) =
h_2(\tau_2,\mathcal{B})
$$
for any $\mathcal{B}$, hence $h_1(\tau_1) \geq h_2(\tau_2)$.

\bigskip

\section{MV-algebras with product}

MV-algebra with product (see [8], [13], [14]) is a pair $(M,.)$, where $M$ is an MV-algebra, and $\cdot$ is a commutative and associative operation on $M$ satisfying the following conditions:

\begin{itemize}
\item[(i)]  $ u.a = a; $
\item[(ii)]
$ if \ a,b,c \in M, a + b \leq u,\text \ {then} \
 c.a + c.b \leq u, \text  \ {and} \ c.(a + b) = c.a + c.b. $
\end{itemize}

In MV-algebras with product a suitable theory of entropy of the Kolmogorov
type has been constructed in [12]. Of course, the construction is different.
The main idea is the following. If $\mathcal{A} = \{a_1,...,a_k\}, \mathcal{B} =
 \{b_1,...,b_l\}$ are two partitions of unity, then one defines
$$
\mathcal{A} \vee \mathcal{B} = \{ a_i.b_j; i=1,...,k, j=1,...,l\}.
$$
The entropy is defined by the formula
$$
\overline h (\mathcal{A},\tau) = \lim_{n \to \infty} \frac 1 n
H(\bigvee_{i=0}^{n-1} \tau^i(\mathcal{A})).
$$
Evidently $\mathcal{A} \vee \mathcal{B}$ is a common
refinement of $\mathcal{A}$ and $\mathcal{B}$. Therefore
$$
H_n(\mathcal{A},\tau) \leq H(\bigvee_{i=0}^{n-1}\tau^i(\mathcal{A})),
$$
and
$$
h(\mathcal{A},\tau) = \lim_{n \to \infty} \frac 1 n H_n(\mathcal{A},\tau) \leq
\lim_{n \to \infty} \frac 1 n H(\bigvee_{i=0}^{n-1}\tau^i(\mathcal{A})) = \overline h(\mathcal{A},\tau).
$$

We want to study some relations between $h(\mathcal{A},\tau)$ and $\overline h (\mathcal{A},\tau)$. Of course,
we shall assume that the given MV-algebra is \it divisible, \rm i.e. to any $ a \in M$ and
any $n \in N$ there exists $b \in M$ such that $ a = n.b$. As it has been shown in [4] any divisible $\sigma$-complete
MV-algebra  $M$ admits a product such that $(M,.)$ is an MV-algebra with product.

\bf Theorem 3. \it Let $M$ be a divisible $\sigma$-complete MV-algebra, m be continuous, po\-si\-tive state on M
(i.e. $m(a) = 0 \Longrightarrow a = 0$). If $\mathcal{A} = \{a_1,...,a_k\}, \mathcal{B} = \{b_1,...,b_l\}$ are partitions of unity,
 and $\mathcal{A}$ consists of idempotents (i.e. $a_i \oplus a_i = a_i$ for all i), then $\mathcal{A} \vee
 \mathcal{B}$ is the unique common refinement of $\mathcal{A}$ and $\mathcal{B}$.

\bigskip

\rm Proof. Our main tool is the Dvure\v censkij and Mundici modification of the Loomis - Sikorski theorem
([2], [10]). By [3] Theorem 7.2.6 there exists a compact Hausdorff space $\Omega$,
a family $\mathcal{T}$ of functions $f:\Omega \rightarrow [0,1]$ and a mapping
$\psi : \mathcal{T} \rightarrow M$
satisfying the following properties:
\begin{itemize}
\item[(i)]  $ \mathcal{T}$ is a tribe, i.e. $1_{\Omega} \in \mathcal{T}, 1 - f \in \mathcal{T}$
whenever $f \in \mathcal{T}$, and $\min\{\sum_{i=1}^\infty, 1\} \in \mathcal{T} $
, whenever $f_n \in \mathcal{T} (n=1,2,...)$;
\item[(ii)]
 the natural product $f.g \in \mathcal{T}$, whenever $f \in \mathcal{T}, g \in \mathcal{T}$;
\item[(iii)] $\psi$ is an MV-$\sigma$-homomorphism from $\mathcal{T}$ onto $M$ preserving the product
in $\mathcal{T}$.
\end{itemize}

Since $M$ is divisible and $\sigma$-complete, the tribe $\mathcal{T}$ contains all constant functions. Indeed, by
[3] theorem 7.1.22 $\Omega$ consists of all functions $m:M \rightarrow [0,1]$ satisfying the identity
$m(a \oplus b) = \min (m(a) + m(b), 1)$, and the tribe $\mathcal{T}$ is generated by the family of all functions
 $\widehat a:\Omega \rightarrow [0,1]$ defined by $\widehat a (m) = m(a), a \in M.$  Of course, divisibility
 of $M$ implies for any $q \in N$ the existence of $\frac 1 q u$ such that $q m(\frac 1 q u) = 1$. If $r \in Q, r = \frac p q$,
 then $\widehat{ru}(m) = m(ru) = p m(\frac 1 q u) = p \frac 1 q m(u) = r$, hence the constant function $r_{\Omega}$ belongs to
  $\mathcal{T}$. Finally $\sigma$-completeness of $M$ implies that $\mathcal{T}$ contains all constant functions $\alpha_{\Omega},
\  \alpha \in [0,1]$.

  Since $\mathcal{T}$ contains all constant fuzzy sets, by the Butnariu - Klement theorem ([1] Prop. 3.3.,
  [15] Theorem 7.1.7) $\mathcal{T}$ consists of all $\mathcal{S}$-measurable functions, where
  $\mathcal{S} = \{ A \subset \Omega ; \chi_A \in \mathcal{T}\}$.

  Define $\mu :\mathcal{T} \rightarrow [0,1]$ by $\mu (f) = m(\psi (f))$. Evidently $\mu$ is a continuous
  state on the $MV$-algebra $\mathcal{T}$. There exists a probability measure $P:\mathcal{S} \rightarrow [0,1]$
  such that $\mu (f) = \int_{\Omega} f dP$ for any $f\in \mathcal{T}$.

  Let $\{c_{ij} ; i=1,...,k, j=1,...,l\}$ be any common refinement of $\mathcal{A}$ and $\mathcal{B}$.
  Since $\psi$ is epimorphism, there are $f_1,...,f_k, g_1,...,g_l, h_{ij} \in \mathcal T$ such
  that $\psi (f_i) = a_i, \psi(g_j) = b_j, \psi(h_{ij}) = c_{ij}.$ Since $a_i$ are idempotent,
  it is not difficult to see ([3] Prop. 7.1.20) that $f_i = \chi_{A_i}$.
  We have P-almost everywhere
$$
  \sum_{j=1}^k h_{ij} = \chi_{A_i}, \sum_{i=1}^k h_{ij} = g_j, i = 1,...,k, j=1,...,l.
$$
  Therefore
$$
  h_{ij} = \chi_{A_i} . g_j,
$$
P-a.e., hence
$$
  c_{ij} = \psi (h_{ij}) = \psi (\chi_{A_i}). \psi (g_j) = a_i.b_j.
$$

\bigskip

\bf Corollary. \it Let the assumptions of Theorem 3 be satisfied.
 If $\mathcal{A}$ consists of idempotent elements, then
 $$
 h(\mathcal{A}, \tau) = \overline h (\mathcal A, \tau).
 $$

 \bf Theorem 4. \it Let the assumptions of Theorem 3 be satisfied. If $\mathcal {A}$
 consists of idempotent elements, then for any partition $\mathcal{B}$ there holds
 the inequality
 $$
 h(\mathcal{B},\tau) \leq h(\mathcal{A},\tau) + H(\mathcal B \Vert \mathcal{A})
 $$
 where
 $$
 H(\mathcal{B}\Vert \mathcal{A}) = \sum_{i=1}^k \sum_{j=1}^l m(b_j) \varphi(\frac {m(a_i.b_j)} {m(b_j)}).
 $$

 \rm Proof. It follows by [12] Prop. 8 and the above Corollary.

 \bigskip

\bf Remark. \rm Theorem 4 presents a basic tool for counting entropy in some structures.
For details see e.g. [15], Chapter 10.

\bigskip

\bf Acknowledgement. \rm  The author thanks to Professors Daniele Mundici
and Anatolij Dvure\v censkij for many fruitful discussions and suggestions.
The paper was supported by Grant VEGA 1/9056/02.

\end{document}